# Modeling of mixed-mechanism stimulation for the enhancement of geothermal reservoirs

Hau Trung Dang │ Eirik Keilegavlen │ Inga Berre

Center for Modeling of Coupled Subsurface Dynamics, Department of Mathematics, University of Bergen

**Correspondence**

Hau Trung Dang, University of Bergen, Postboks 7803, 5020 Bergen, Norway

Email:

Hau.Dang@uib.no

dtrhau@gmail.com

**Abstract**

Hydraulic stimulation is a critical process for increasing the permeability of fractured geothermal reservoirs. This technique relies on coupled hydromechanical processes induced by reservoir stimulation through pressurized fluid injection into the rock formation. The injection of fluids causes poromechanical stress changes that can lead to the dilation of fractures due to fracture slip and to tensile fracture opening and propagation, so-called mixed-mechanism stimulation. The effective permeability of the rock is particularly enhanced when new fractures connect with pre-existing fractures. Mixed-mechanism stimulation can significantly improve the productivity of geothermal reservoirs, and the technique is especially important in reservoirs where the natural permeability of the rock is insufficient to allow for commercial flow rates.

This paper presents a modeling approach for simulating the deformation and expansion of fracture networks in porous media under the influence of anisotropic stress and fluid injection. It utilizes a coupled hydromechanical model for poroelastic, fractured media. Fractures are governed by contact mechanics and allowed to grow and connect through a fracture propagation model. To conduct numerical simulations, we employ a two-level approach, combining a finite volume method for poroelasticity with a finite element method for fracture propagation. The study investigates the impact of injection rate, matrix permeability, and stress anisotropy on stimulation outcomes. By analyzing these factors, we can better understand the behavior of fractured geothermal reservoirs under mixed-mechanism stimulation.

## 1. Introduction

Hydraulic stimulation plays a critical role in facilitating the production of geothermal energy in low-permeability igneous rocks. Its main goal is to increase reservoir permeability to achieve flow rates that are economically feasible for commercial production.[1–3] Hydraulic stimulation can be performed at different fluid pressures. High pressures exceeding the minimum principal stress are used to propagate hydraulic fractures, while elevated but lower pressures can cause hydro-shearing of pre-existing fractures as their frictional resistance to slip is exceeded.

In conventional hydraulic fracturing, mixtures of liquid and small insoluble particles are injected at pressures exceeding the tensile strength of the rock to increase reservoir permeability. A high-pressure injection may cause stress concentration at the fracture tip that can trigger tensile fracture propagation.[4] Propagating fractures may connect with pre-existing fractures[5] and thereby increase the fluid flow. When the hydraulic pressurization is reduced, small insoluble particles are retained in the opening of the fracture and, hence, maintain increased permeability. When applied to geothermal reservoirs, this process risks thermal short-circuiting and corresponding low temperatures of the produced fluid.[6,7]

Injections at pressures below the minimum principal stress have been shown to be an efficient mechanism for stimulating larger volumes of rock if the reservoir is characterized by pre-existing fractures and faults and high-stress anisotropy. In this case, poromechanical stress changes induced by fluid injection can cause fracture slip and corresponding shear dilation due to the sliding of rough fracture surfaces against each other. Shear dilation can strongly enhance fracture permeability.[8,9] For injections at pressures close to and above the minimum principal stress, the deformation of pre-existing fractures combines with the propagation of wing cracks toward the direction of maximum principal stress.[3,10–13] When a propagating fracture reaches another pre-existing fracture, there are no pressure concentration and low tensile stress at the tip; thus, propagation is arrested.[5] The pressure increase due to injection can then extend to the newly connected fracture, potentially causing shear slip or tensile opening and the formation of new wing cracks. As a result, the development of complex fracture networks created by connecting newly formed wing cracks to pre-existing fractures enhances the permeability of the geothermal reservoir. This mechanism of hydraulic stimulation, combining shear-dilation and propagation of fractures, is referred to as mixed-mechanism stimulation.[13–15] However, the complex dynamics of stress redistribution related to mixed-mechanism stimulation and how it interacts with pre-existing fractures are not well understood.

Numerical modeling can be employed to study the interaction between fluid flow through fractured rock and the poromechanical deformation of the rock, including fracture deformation and propagation. The complexity of the coupled processes makes it difficult to include all such effects, and thus it is common to apply simplified models that consider only a subset of the processes. For instance, modeling of tensile fracturing of poroelastic media caused by high injection pressure while neglecting the effects of shear slip, contact, and friction has been widely reported.[16–19] Several studies have further investigated the extension of pre-existing fracture networks in porous media resulting from fluid injection. However, these studies have either neglected friction and contact mechanics at fracture interfaces[20–22] or forced fractures to propagate along predefined paths.[11,21]

Recently, the authors proposed a new methodology to simulate fluid flow, matrix deformation, fracture slip, and fracture propagation in porous media as a result of fluid injection.[10] Specifically, a mathematical model was developed based on the mixed-dimensional discrete fracture matrix (md-DFM) conceptual model that combined the explicit representation of major fractures with a continuum representation of the surrounding medium. This model utilized a co-dimension-one representation of the fractures. Hence, for a two-dimensional (2D) domain, fractures were represented as one-dimensional (1D) lines, with a longitudinal parameter representing fracture apertures. The model allowed for the application of fracture contact mechanics, including frictional sliding and shear-dilation of fractures and tensile fracture opening. The framework was designed as a two-level method, with local computation of fracture propagation around individual tips split separate from global computations of flow and poromechanical deformation of the fractured rock. The coupling strength between the local and global models was a user-controlled parameter that allowed users to balance simulation accuracy and computational cost.

This study uses the approach proposed by Hau et al.[10] to further investigate the mixed-mechanism stimulation of fractured rock under anisotropic stresses. It explores how fluid injection can change the effective poroelastic stress regime, resulting in fracture slip and dilation as well as tensile fracture propagation. The study examines how stimulation outcomes are affected by the injection rate, matrix permeability, and stress anisotropy. Specifically, the study considers fracture coalescence, which creates new, dominant flow paths.

The paper is organized as follows. Section 2 presents the mathematical model for mixed-mechanism stimulation of a fractured geothermal reservoir. In Section 3, we describe the numerical approach used to simulate the behavior of the reservoir under stimulation. Section 4 presents the results of several numerical test cases, which provide insights into the role of mixed-mechanism stimulation in enhancing reservoir permeability. Finally, in Section 5, we present our conclusions and provide remarks about the implications of our findings.

## 2. Mathematical model

This section presents the governing equations that model fluid flow and deformation in fractured porous media. Additionally, we introduce a mathematical model for fracture contact mechanics, propagation, and coalescence. These equations are essential for developing a simulation model that accurately captures the behavior of fractured geothermal reservoirs under mixed-mechanism stimulation. By modeling the coupling of fluid flow, rock deformation, and fracture growth, we can better understand the impact of stimulation on the reservoir. The numerical simulations described in later sections of the paper are based on the mathematical models presented in this section.

### 2.1. Fluid flow and poroelastic deformation of the matrix and fracture

The md-DFM conceptual model for a 2D fractured porous media domain was employed in this paper. By using the md-DFM model, we divide the domain into three subdomains: a 2D host medium denoted by $\Omega^M$, a set of fractures represented as 1D objects and denoted by $\Omega^F$, and fracture intersections represented as points and denoted by $\Omega^I$. The boundaries of $\Omega^M$ and $\Omega^F$ are denoted by $\partial\Omega^M$ and $\partial\Omega^F$, respectively, while $\Gamma$ represents the interfaces between the host medium and fractures. When necessary, to denote the interfaces at the different sides of a fracture, we use superscripts $\pm$ on $\Gamma$. The interfaces between $\Omega^F$ and $\Omega^I$ are denoted by $\Lambda$, where the superscript $i$ is used on $\Lambda$ when necessary to denote the interface between $\Omega^I$ and a specific fracture indexed by $i$. Figure 1 provides an illustration of the model.

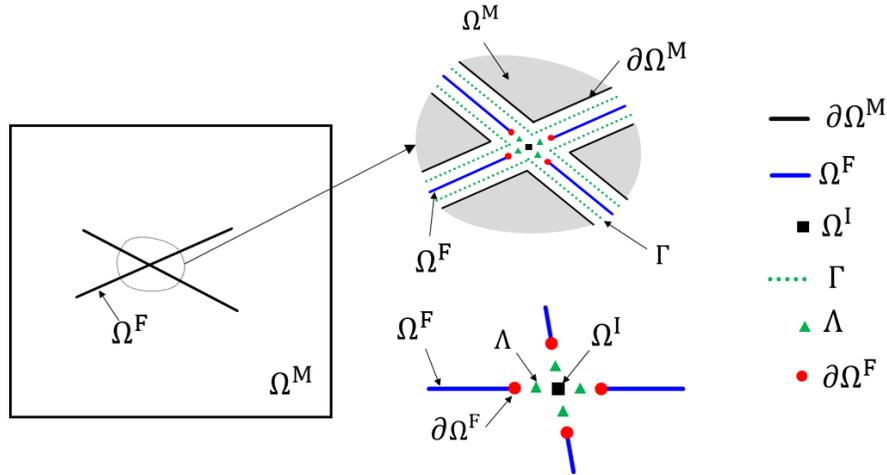

Figure 1. Illustration of a host medium $\Omega^M$, fractures $\Omega^F$, intersection $\Omega^I$, and interfaces between higher- and lower-dimensional domains, denoted by $\Gamma$ and $\Lambda$, respectively. In the detailed images to the right of the general figure on the left, the different domains and interfaces are separated for illustration purposes.

To facilitate coupling between the subdomains, projection operators $\Pi_{[-]}^{[-]}$ are introduced.[23] The illustration of these operators is given in Figure 2, where the subscripts of $\Pi$ indicate the origin, while the superscripts indicate the destination of the projection.

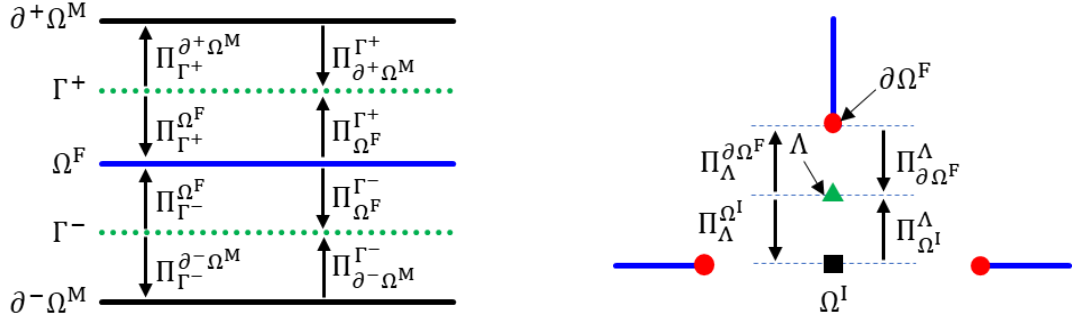

a) Projection operators between $\Omega^M$ and $\Omega^F$     b) Projection operators between $\Omega^F$ and $\Omega^I$

Figure 2. Illustration of projection operators between subdomains.

In our model, we assume that the porous media domain is deformable and that its mechanical properties are elastic, isotropic, and homogeneous. We assume that the fluid is a single phase and slightly compressible. The permeability is allowed to be heterogeneous. The governing equations can be given as follows:

$$\nabla \cdot \boldsymbol{\sigma} = \mathbf{b}, \qquad\qquad \text{on } \Omega^M \qquad (1)$$

$$\boldsymbol{\sigma} = \mathbf{C}\nabla\mathbf{u} - \alpha p \mathbf{I}, \qquad\qquad \text{on } \Omega^M \qquad (2)$$

$$\alpha \frac{\partial(\nabla \cdot \mathbf{u})}{\partial t} + \left(\phi c_p + \frac{\alpha - \phi}{K}\right)\frac{\partial p}{\partial t} + \nabla \cdot \mathbf{q} = q^0, \qquad\qquad \text{on } \Omega^M \qquad (3)$$

$$\mathbf{q} = -\frac{1}{\mu}\begin{bmatrix} \kappa_{xx} & 0 \\ 0 & \kappa_{yy} \end{bmatrix}\nabla p, \qquad\qquad \text{on } \Omega^M \qquad (4)$$

$$\frac{\partial a}{\partial t} + a c_p \frac{\partial p_F}{\partial t} + \nabla_\parallel \cdot \mathbf{q}_F - \Pi_{\Gamma^+}^{\Omega^F}\lambda^+ - \Pi_{\Gamma^-}^{\Omega^F}\lambda^- = q_F^0, \qquad\qquad \text{on } \Omega^F \qquad (5)$$

$$\mathbf{q}_F = -\frac{\kappa a}{\mu}\nabla_\parallel p_F, \qquad\qquad \text{on } \Omega^F \qquad (6)$$

$$\frac{\partial(a_I^2)}{\partial t} + a_I{}^2 c_p \frac{\partial p_I}{\partial t} - \sum_{i=1}^{N}\Pi_{\Lambda^I}^{\Omega^I}\eta_i = q_I^0, \qquad\qquad \text{on } \Omega^I \qquad (7)$$

where $\mathbf{u}, \boldsymbol{\sigma}, p$, and $\mathbf{q}$ denote displacements, stress, pore pressure, and flux on $\Omega^M$. The source terms for the mass conservation in the subdomains for the matrix, fractures, and fracture intersections are denoted by $q^0$, $q_F^0$, and $q_I^0$, respectively. The flux and pressure in the fracture subdomains are denoted by $\mathbf{q}_F$ and $p_F$, respectively. The terms $\lambda^\pm$ are variables that represent the flux from the matrix to the fracture at each side of the fracture. The aperture of the fracture is $a$, and for the fracture intersection, $\Omega^I$, the aperture, $a_I$, is taken to be the average of the apertures of the intersecting fractures. The term $\eta_i$ is a variable that represents the flux from fracture $i$ to $\Omega^I$, and $N$ is the number of intersecting fractures around $\Omega^{\mathrm{I}}$. The fracture aperture is a function determined based on the residual aperture and normal displacement jump, such that:

$$a = a_0 + [\![\mathbf{u}]\!]_{\mathrm{n}}, \qquad\qquad\qquad \text{on } \Omega^{\mathrm{F}} \quad (8)$$

where $a_0$ denotes the residual aperture in the undeformed state, and $[\![\mathbf{u}]\!]_{\mathrm{n}}$ represents the displacement jump in the normal direction over $\Omega^{\mathrm{F}}$, in which the displacement jump is defined by:

$$[\![\mathbf{u}]\!] = \mathbf{u}|_{\Gamma^-} - \mathbf{u}|_{\Gamma^+}, \qquad\qquad\qquad \text{on } \Omega^{\mathrm{F}} \quad (9)$$

where $\Gamma$ is the interface between $\Omega^{\mathrm{M}}$ and $\Omega^{\mathrm{F}}$. The other parameters in the above equations are given in Table 1.

Table 1. The parameters used in the governing equations.

| Notation | Description | Notation | Description |
|---|---|---|---|
| $\mathbf{C}$ | stiffness matrix | $c_p$ | fluid compressibility |
| $\phi$ | matrix porosity | $\mu$ | fluid viscosity |
| $\kappa_{xx}, \kappa_{yy}$ | permeability of the porous matrix | $N$ | number of intersecting fractures |
| $\zeta$ | inflow from the matrix to the fracture | $K$ | bulk modulus |
| $\kappa$ | fracture permeability | $\mathbf{b}$ | body forces around $\Omega^{\mathrm{I}}$ |
| $\nabla, \nabla_{\parallel}$ | gradient operators | tr | trace operator |

To fully represent the physical system, it is necessary to incorporate the coupling between subdomains into the mathematical model. First, the coupling between $\Omega^M$ and $\Omega^F$ is defined by:

$$\mathbf{q} \cdot \mathbf{n}|_{\partial^{\pm}\Omega^{\mathrm{M}}} = \Pi_{\Gamma^{\pm}}^{\partial^{\pm}\Omega^{\mathrm{M}}} \lambda^{\pm}, \qquad\qquad \text{on } \partial\Omega^{\mathrm{M}} \quad (10)$$

$$\lambda^{\pm} = -\frac{\kappa}{\mu}\left(\frac{\Pi_{\Omega^{\mathrm{F}}}^{\Gamma^{\pm}} p_F - \Pi_{\partial^{\pm}\Omega^{\mathrm{M}}}^{\Gamma^{\pm}} \mathrm{tr}^{\pm} p}{a/2}\right), \qquad\qquad \text{on } \Gamma^{\pm} \quad (11)$$

where Eq. (10) indicates the balance of flux between the matrix and fracture. The coupling between $\Omega^{\mathrm{F}}$ and $\Omega^{\mathrm{I}}$ is given by:

$$\mathbf{q}_f \cdot \mathbf{n}|_{\partial\Omega_i^{\mathrm{F}}} = \Pi_{\Lambda^{\mathrm{i}}}^{\Omega^{\mathrm{F}}} \eta_i, \qquad\qquad \text{on } \partial\Omega_i^{\mathrm{F}} \quad (12)$$

$$\eta_i = -\frac{\kappa a_I}{\mu}\left(\frac{\Pi_{\Omega^{\mathrm{I}}}^{\Lambda^{\mathrm{i}}} p_I - \Pi_{\Omega^{\mathrm{F}}}^{\Lambda^{\mathrm{i}}} p_F}{a_I/2}\right). \qquad\qquad \text{on } \Lambda_i \quad (13)$$

The governing equations presented here are comprehensive, as they describe the mechanisms operating in each subdomain and consider their interactions.

## 2.2. Fracture contact mechanics

In the context of hydromechanical coupled processes, fractures are assumed to be in one of three states: closed and sticking (with no shear displacement), closed and slipping, or open. The interactions between the fracture surfaces are governed by fracture contact mechanics. In the following, the fracture contact mechanics model is considered independently in the normal and tangential directions. First, the normal opening of the fracture is governed by a non-penetration condition written in Karush-Kuhn-Tucker (KKT) form[24] as:

$$[\![\mathbf{u}]\!]_n - g \geq 0, \quad f_n \leq 0, \quad ([\![\mathbf{u}]\!]_n - g)f_n = 0. \qquad \text{on } \Omega^F \quad (14)$$

Here, $f_n$ represents the contact traction in the normal direction, and $g$ is a gap function defined by:

$$g = -\tan(\psi)\|[\![\mathbf{u}]\!]_\tau\|, \qquad \text{on } \Omega^F \quad (15)$$

where $\psi$ is the dilation angle and $[\![\mathbf{u}]\!]_\tau$ is the displacement jump in the tangential direction. The gap function in Eq. (15) accounts for the dilation of the fracture resulting from tangential slip while maintaining contact between the fracture surfaces. This feature enables the enhancement of permeability in the fracture due to shear dilation.

The tangential motion of the fracture is modeled as a frictional contact problem given by:

$$\begin{aligned} &|f_\tau| \leq -\mu_s f_n, \\ &|f_\tau| < -\mu_s f_n \rightarrow [\![\dot{\mathbf{u}}]\!]_\tau = 0, \\ &|f_\tau| = -\mu_s f_n \rightarrow \exists \varepsilon \in \mathbb{R}, \; f_\tau = -\varepsilon[\![\dot{\mathbf{u}}]\!]_\tau, \end{aligned} \qquad \text{on } \Omega^F \quad (16)$$

where $\mu_s$ represents the friction coefficient and $\dot{\mathbf{u}}$ is the derivative of $\mathbf{u}$ with respect to time. The contact traction in the tangential direction, $f_\tau$, contains directional information, and is therefore a vector despite the fracture being 1D.

Traction on the fracture surfaces balances the pressure in the fracture by Newton's third law and can be expressed as

$$\mathbf{f}_+ = \left(\Pi_{\partial^+ \Omega^M}^{\Gamma^+} \boldsymbol{\sigma} + \mathbf{I}\, \alpha_f\, \Pi_{\Omega^F}^{\Gamma^+} p_F\right), \qquad \text{on } \Gamma^+ \quad (17)$$

$$\mathbf{f}_- = -\left(\Pi_{\partial^- \Omega^M}^{\Gamma^-} \boldsymbol{\sigma} + \mathbf{I}\, \alpha_f\, \Pi_{\Omega^F}^{\Gamma^-} p_F\right). \qquad \text{on } \Gamma^- \quad (18)$$

The tractions on $\Gamma^\pm$ are related to the contact traction vector $\mathbf{f} = (\boldsymbol{f}_\tau, f_n)$ by $\mathbf{f}_\pm = \pm \Pi_{\Omega^F}^{\Gamma^\pm}(\mathbf{R}\mathbf{f})$, where $\mathbf{R}$ is a rotation matrix from the local $(\tau, n)$ to the global $(x, y)$ coordinate system. Eqs. (17) and (18) indicate that the traction on the fracture surfaces is caused not only by the matrix deformation and pressure but also by pressure in the fracture.

## 2.3. Fracture propagation

We combine the maximum tangential stress criterion[25] and Paris's law[26] to determine the onset of fracture propagation as well as the propagation direction and length. The maximum tangential stress criterion assumes that a fracture propagates when the maximum tangential stress in the process zone around a fracture tip exceeds a critical value defined as:

$$K_{\mathrm{I}} \cos^3 \frac{\theta}{2} - \frac{3}{2} K_{\mathrm{II}} \cos \frac{\theta}{2} \sin \theta \geq K_{\mathrm{IC}}, \tag{19}$$

The direction of propagation is that of the maximum tangential stress given by:

$$\theta = 2 \tan^{-1} \left( \frac{K_{\mathrm{I}}}{4 K_{\mathrm{II}}} \pm \frac{1}{4} \sqrt{\left( \frac{K_{\mathrm{I}}}{K_{\mathrm{II}}} \right)^2 + 8} \right), \tag{20}$$

$$K_{\mathrm{II}} \left( \sin \frac{\theta}{2} + 9 \sin \frac{3\theta}{2} \right) < K_{\mathrm{I}} \left( \cos \frac{\theta}{2} + 3 \cos \frac{3\theta}{2} \right), \tag{21}$$

where $K_{\mathrm{I}}$ and $K_{\mathrm{II}}$ are the stress intensity factors (SIFs). If more than one crack grows simultaneously, then the tips in the fracture with higher energy advance farther than the others, with a distribution given by the Paris-type law,[26]

$$l_{\mathrm{adv}}^i = l_{\max} \left( \frac{G_i}{\max(G_i)} \right)^{0.35}, \tag{22}$$

where $l_{\mathrm{adv}}^i$ and $G_i$ are the propagation length and energy release for tip $i$, respectively.[27] By Eq. (22), the increment for each tip is limited by a preset value, $l_{\max}$.

A propagating fracture may reach and coalesce with another fracture in a T-type connection. This leads to the formation of a new intersection point that is added to $\Omega_{\mathrm{I}}$ and new connections between the merged fractures and the intersection.

## 3. Discretization method

In this section, we describe a numerical approach for discretizing the mathematical model presented in Section 2. As the model depends on both space and time variables, both variables must be discretized. Since the mathematical model contains only the first derivative with respect to time, time discretization can be achieved using the backward Euler method. However, the model is more complex regarding spatial variations, which can be dealt with by the two-level simulation recently proposed by Hau et al.[10]

The motivation for using the two-level simulation approach is to balance computational cost and simulation accuracy. Specifically, poroelastic deformations with fracture contact mechanics, but without fracture propagation, are assumed to be quasi-static and are treated using a relatively coarse grid. In contrast, a locally refined grid around the fracture tip is needed to accurately capture the interaction between fracture propagation and local stress variations. If a fracture propagates and exceeds a certain threshold length, then the geometry of the fracture network and the solution are updated in the coarse-level domain for the next time step. A brief description of this approach is provided below; for more information, we refer to Hau et al.[10]

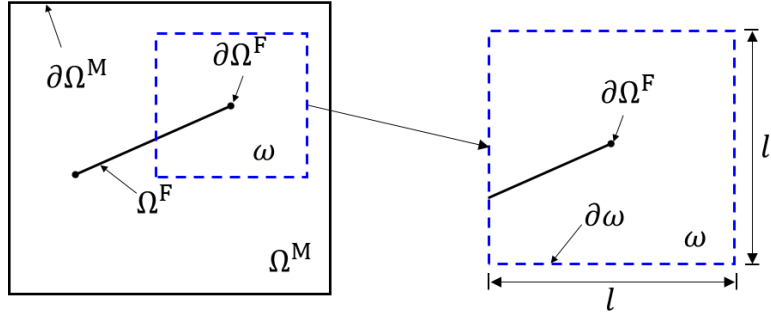

Figure 3. Illustration of a fracture, $\Omega^{\mathrm{F}}$, and a fine-level domain $\omega$, adapted from Hau et al. (2022).[10]

The computational domain is divided into a coarse-level domain that matches the entire domain and smaller fine-level domains with size $l$ that surround the fracture tips. The coarse-level and fine-level domains are denoted by $\Omega^{\mathrm{M}}$ and $\omega$, respectively, as illustrated in Figure 3. These domains are discretized using triangular cells with grid sizes $\Delta H$ and $\Delta h$ for $\Omega^{\mathrm{M}}$ and $\omega$, respectively. The grids conform to the fractures so that fractures coincide with grid faces, and nodes and faces are split along the fractures. To best represent fracture paths in the grids and avoid excessive computational cost while ensuring an accurate numerical solution at the relevant scale of the model, an adaptive remeshing technique is employed.[27] This technique uses finer cells around fracture tips in both coarse-level and fine-level grids to sufficiently capture the details of fracture propagation. Additionally, to ensure the stability of the simulation, the resolution of the fine-level grid is set to be finer than that of the coarse-level grid, i.e., $\Delta h = \varepsilon_m \Delta H$ with $\varepsilon_m \leq 1$.

When none of the fine-level domains intersect with neighboring fractures, the coarse-level and fine-level domains are defined differently. However, for technical reasons, our implementation cannot handle fine-scale domains that contain multiple fractures. Therefore, when there is an intersection between a fine-level domain and neighboring fractures, the fine-level domain is defined to be identical to the coarse-level domain. Nonetheless, we emphasize that the proposed approach is still applicable for much larger domains than those demonstrated in this paper.

## 3.1. Two-level discretization

The poroelastic deformation model presented in subsections 2.1 and 2.2 is discretized based on the coarse-level grid. Specifically, the governing equations in subsection 2.1 are discretized using a finite volume approach with a multi-point flux approximation and a multi-point stress approximation,[28,29] while the fracture contact mechanics presented in subsection 2.2 are discretized by an active set method.[24,30,31] The solution at this level provides the deformation and fluid pressure in the poroelastic domain and determines fracture mechanical behavior, whether the fracture is open, closed and sticking or closed and slipping.

The fine-level domain is responsible for evaluating fracture propagation at each time step. To do this, we combine Eqs. (1) and (2) and assume that the fine-level domain behaves similarly to a linearly elastic medium governed by:

$$\nabla \cdot (\mathbf{c}\nabla_s \mathbf{u}_l) + \mathbf{b} = 0, \tag{23}$$

where $\mathbf{u}_l$ is the deformation in the fine-level domain, $\mathbf{b} = -\nabla \cdot (\alpha p \mathbf{I})$ is the body force caused by pressure from the coarse-level domain, and $\mathbf{c}$ is the stiffness tensor. The boundary conditions for the fine-level problem, i.e., defined at $\partial \omega$, are set according to the coarse-level state. To solve Eq. (23), we use a combination of the $\mathcal{P}_2$ finite element method and quarter point elements to accommodate the stress singularity at the fracture tip.[32,33] The solution obtained is then used to compute stress intensity factors (SIFs) and determine whether a fracture will propagate and, if so, where and how far it will go, as described in subsection 2.3. The maximum increment of fracture is set to the fine-level grid size, i.e., $l_{\max} = \Delta h$.

### 3.2. Coupling between coarse-level and fine-level solutions

To establish the numerical coupling between the coarse-level and fine-level domains, it is necessary to project the displacements from the coarse-level to the fine-level domain boundaries and compress the fine-level updates to the fracture geometry in the coarse-level grid. These projections can be achieved using three mapping processes: cell center to cell center (C2C), node to node (N2N), and cell center to node (C2N).[10] Additionally, updating the coarse-level fracture path is necessary if the propagation in the fine-level domain is sufficiently significant to cause a considerable change in the coarse-level grid. To accomplish this, we denote $|\Delta \omega^F|$ as the total propagation length in a fine-level domain. If $|\Delta \omega^F|$ exceeds $\varepsilon_p \Delta H$, with $\varepsilon_p$ being a propagation factor, the coarse-level fracture is extended using a linear approximation of $\Delta \omega^F$, and the coarse-level grid is updated.

### 3.3. Fracture coalescence

This paper models the fracture intersection by a T-type connection. As illustrated in Figure 4 (a), when the distance between a propagating crack tip and a boundary or another fracture is less than the grid size around the tip, the two fractures are assumed to be connected. A connection point is identified by projecting the fracture tip onto the boundary, resulting in point A. Point B is then defined as the projection of point A to the opposite side of the connected fracture boundary. Finally, the tip of the propagating fracture is split at point A to create a T-type connection, as depicted in Figure 4 (b).

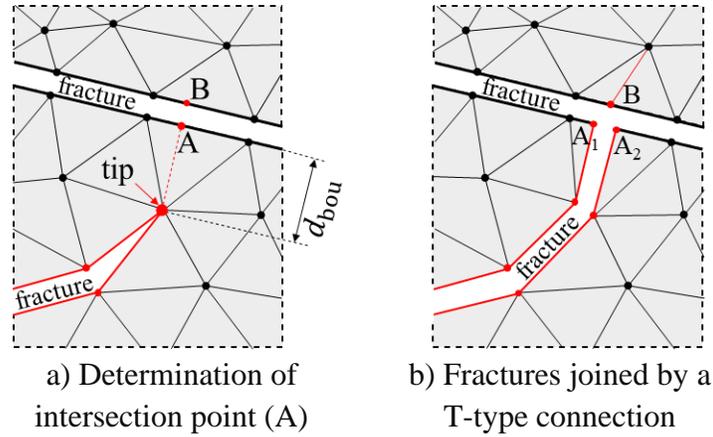

a) Determination of intersection point (A)   b) Fractures joined by a T-type connection

Figure 4. The T-type intersection between fractures or a fracture and boundary. The fracture is widened for illustration purposes.

## 4. Results

The accuracy of simulations of fracture propagation and fluid flow in the fractured porous media domain were verified in previous studies.[10,23,27] The numerical examples in this section aim to show the ability of the proposed model to simulate complex problems, such as multiple fractures deforming, propagating, and connecting in a medium with heterogeneous permeability.

This section presents four numerical examples to investigate the effects of the fluid injection rate, principal stress, permeability, and fracture network on mixed-mechanism stimulation for a fractured low-permeability porous medium representative of an idealized configuration in a geothermal reservoir. Given the limitation of our resources, a relatively small domain with several pre-existing fractures is considered. For all cases, the coordinates of the tips, the material, and the simulation parameters are given in Table 2, Table 3, and Table 4, respectively.

Table 2. Tips coordinates (units: m)

| Tip | $x$ | $y$ | Tip | $x$ | $y$ |
|-----|-----|-----|-----|-----|-----|
| A | 1.00 | 1.15 | B | 1.00 | 0.85 |
| C | 0.85 | 0.97 | D | 1.15 | 1.03 |
| E | 0.65 | 1.10 | F | 0.65 | 0.90 |
| G | 1.40 | 1.06 | H | 1.28 | 0.94 |

Table 3. Material properties

| Parameter | Definition | Value |
|-----------|-----------|-------|
| $E$ | Young's modulus | 40.0 GPa |
| $\nu$ | Poisson's ratio | 0.2 |
| $K_{\mathrm{IC}}$ | fracture toughness | 1.0 MPa $\cdot$ m$^{1/2}$ |
| $\alpha$ | Biot's coefficient in the matrix | 0.8 |
| $\phi$ | material porosity | 0.01 |
| $c_p$ | fluid compressibility | $4.4 \times 10^{-10} \mathrm{Pa}^{-1}$ |
| $\mu$ | viscosity | $1.0 \times 10^{-4} \mathrm{Pa} \cdot \mathrm{s}$ |
| $\mu_s$ | friction coefficient | 0.5 |
| $\psi$ | dilation angle | $1.0^{\circ}$ |
| $a^0$ | initial aperture | 1.0 mm |

Table 4. Simulation parameters

| Parameter | Definition | Value |
|---|---|---|
| $L_x = L_y$ | coarse-level domain size | 2.0 m |
| $l$ | fine-level domain size | 0.1 m |
| $\Delta H$ | coarse-level grid size | 0.02 m |
| $\Delta h$ | fine-level grid size | 0.01 m |
| $\varepsilon_{\mathrm{m}}$ | ratio between coarse-grid and fine-grid sizes | 0.5 |
| $\varepsilon_p$ | propagation factor | 0.5 |
| $\Delta t$ | time step | 0.5 minutes |

## 4.1. Effect of principal stress direction

First, the effect of the principal stress on fracture propagation is investigated. As illustrated in Figure 5, we consider a 2D domain containing two intersecting fractures and the boundary conditions prescribed in this figure. We assume that the matrix permeability of the domain is isotropic and homogeneous, given by $\kappa_{xx} = \kappa_{yy} = 5.0 \times 10^{-20} \mathrm{m}^2$. The fractured porous medium is subject to a stress state imposed orthogonally to the domain. Fluid is injected into the vertical fracture continuously at a constant rate of $Q_0 = 1 \times 10^{-7} \mathrm{m}^2/\mathrm{s}$. Two stress scenarios are considered. For case 1, $\sigma_1 = 2\sigma_2 = 20$ MPa, and for case 2, $2\sigma_1 = \sigma_2 = 20$ MPa. The propagation of the fractures, presented by solid lines, and the fluid flow, described by color, are shown in Figure 6.

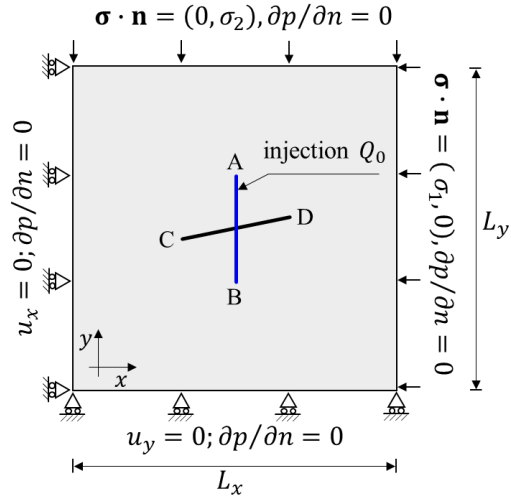

Figure 5. The geometry of model 1.

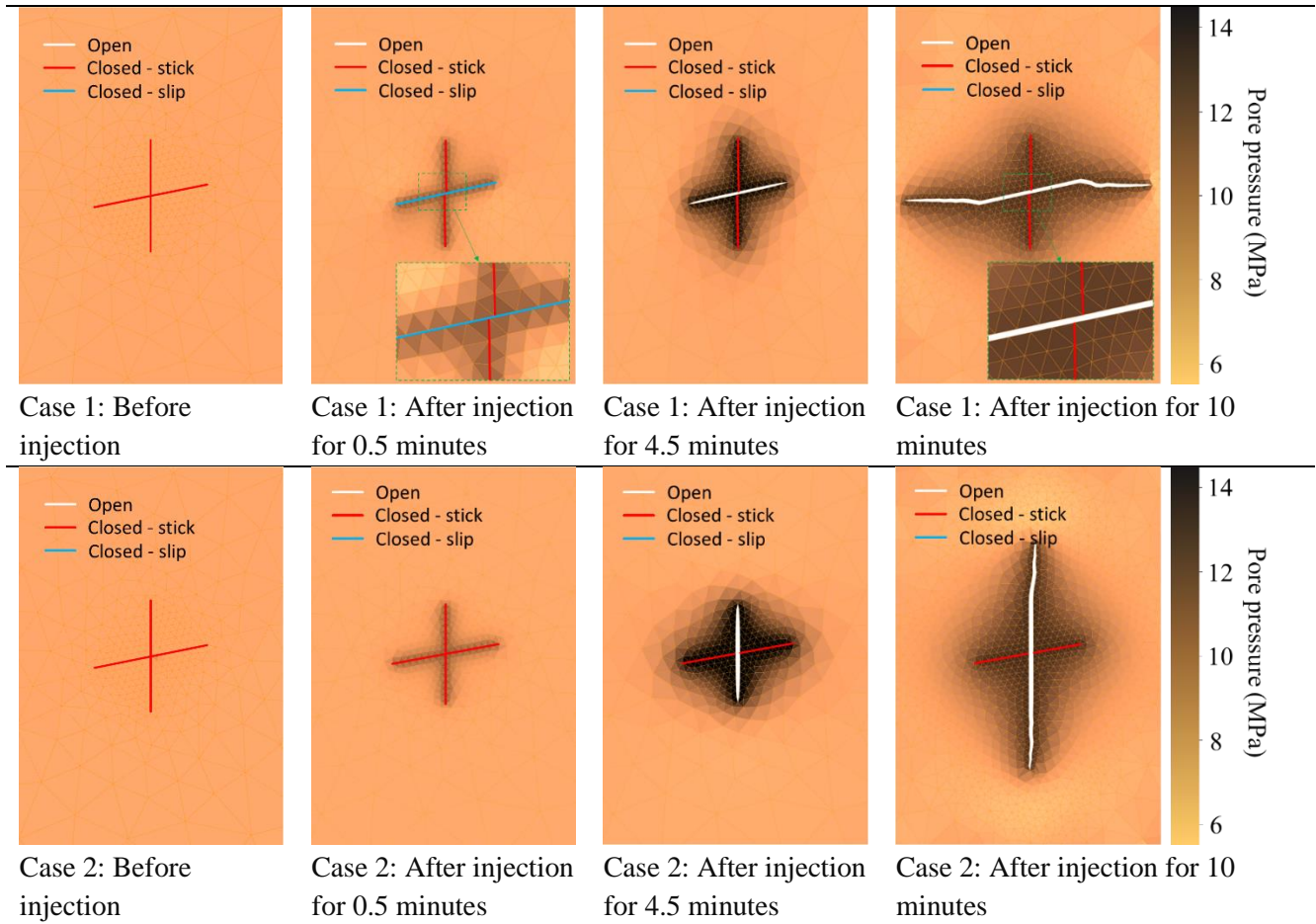

Figure 6. Fracture propagation and pressure evolution in a 2D porous media during fluid injection at rate $Q_0 = 1 \times 10^{-7} \; m^2/s$ into a pre-existing fracture. The solid white lines indicate opening fractures, while the solid red lines indicate closed fractures. The color bar represents pore pressure in MPa.

In both scenarios, pre-existing fractures are closed before fluid is injected due to compressive stress and friction at the fracture interfaces. Depending on the stress regime, the injection can lead to slip in pre-existing fractures. After 0.5 minutes of injection, in case 1, the fracture, which is nearly parallel to the direction of maximum stress, slips. At the same time, for case 2, both fractures remain undeformed, i.e., in the stick mode. In both cases, the vertical fracture is closed and remains in stick mode under compressive stress.

It is well known that fractures propagate toward the direction of maximum principal stress. In case 1, the low injection rate of the fluid does not provide sufficient pressure to induce tensile propagation of the vertical fracture where fluid is injected. However, it does cause shear slip and dilation of the nearly horizontal crossing fracture early in the stimulation process. Continued injection results in wing cracks that appear after 4.5 minutes and propagate in the direction of the maximum principal stress. Thus, this test case demonstrates an example of mixed-mechanism stimulation, where both shear-slip and tensile fracture propagation occur during the stimulation. In case 2, continued fluid injection combined with the shifted stress anisotropy causes the vertical fracture in which the fluid is injected to open. Shear slip does not occur in this case, and tensile propagation of the vertical fracture initiates after 7 minutes of injection once the fluid pressure has built up sufficiently. The simulation also displays the state of fractures, whether they are closed in stick mode, closed in slip mode, or open. A red line indicates a section of a fracture in stick mode, while a light-blue line indicates a section in slip mode. A section of a fracture in open mode is indicated by a solid white line.

## 4.2. Effect of matrix permeability

This study examines the influence of matrix permeability on fracture propagation within a 2D domain. Two distinct permeability regions are investigated, as illustrated in Figure 7. Region 1 is bounded by the curves $c_1: x - (y-1)^2 - 1.2 = 0$, $c_2: x - (y-1)^2 - 1.4 = 0$, and the right boundary, while region 2 is the remainder. The permeability in region 2 is homogeneous and isotropic with values of $\kappa_{xx} = \kappa_{yy} = 5 \times 10^{-20}$ m$^2$. Two simulation cases are conducted, depending on the permeability of region 1. For case 1, $\kappa_{xx} = 5 \times 10^{-20}$ m$^2$ and $\kappa_{yy} = 5 \times 10^{-18}$ m$^2$, while for case 2, $\kappa_{xx} = 5 \times 10^{-20}$ m$^2$ and $\kappa_{yy} = 5 \times 10^{-19}$ m$^2$. Additional parameters used for the simulations are $\sigma_1 = 2\sigma_2 = 20$ MPa and $Q_0 = 1 \times 10^{-7}$ m$^2$/s. The propagation of the fractures and the fluid flow are shown in Figure 8.

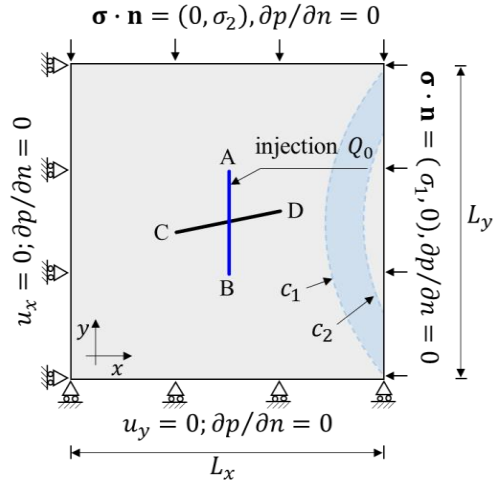

Figure 7. The geometry of model 2.

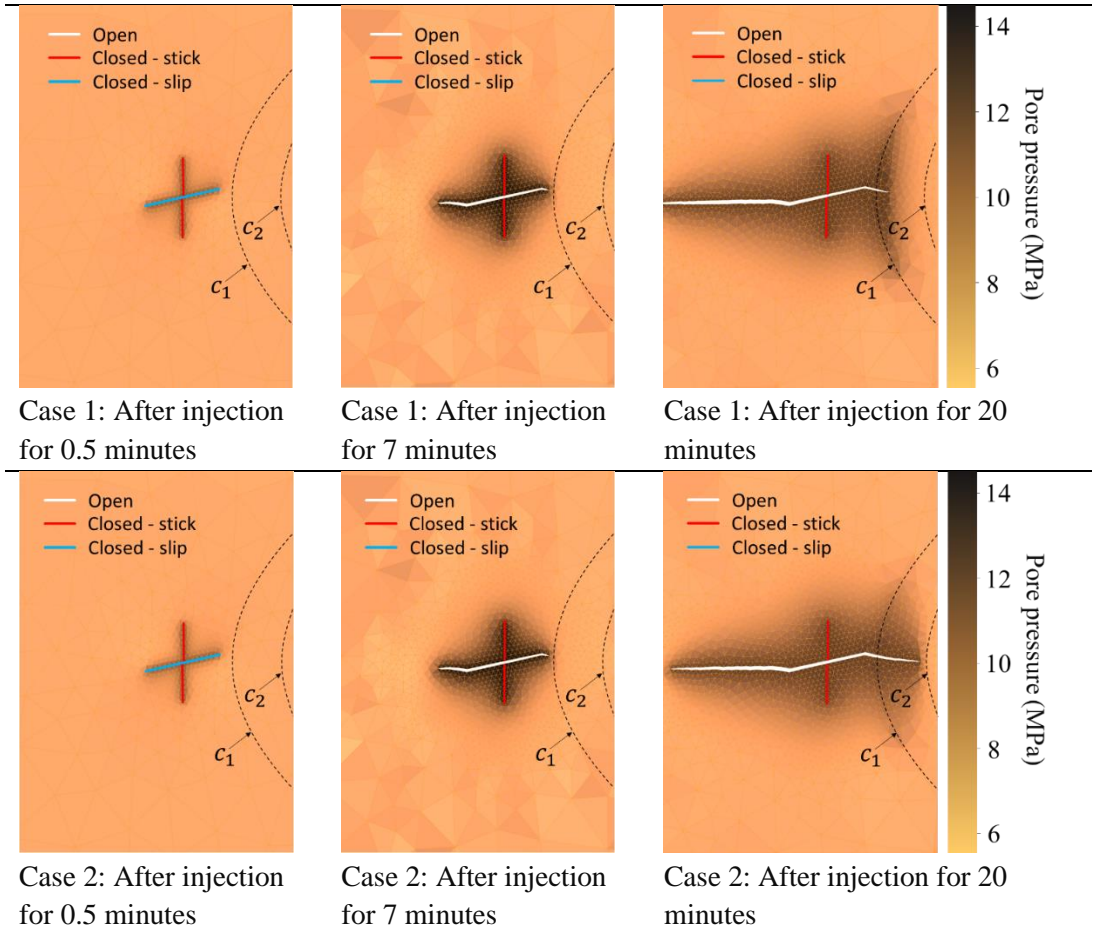

Figure 8. Fracture propagation and pressure evolution in a 2D porous medium during fluid injection, $Q_0 = 1 \times 10^{-7} \ m^2/s$, into a pre-existing fracture. The solid white lines indicate open fractures, while the solid red lines indicate closed fractures. The color bar represents pore pressure in MPa.

The presence of a highly permeable area inhibits fracture growth by preventing fluid pressure from building sufficiently due to fluid leakage into the matrix. Similar to case 1 in example 4.1, the principal stress scenario and fluid injection induce horizontal fracture slip and trigger the appearance of wing cracks after 5.5 minutes of injection. The wing cracks then propagate to opposite sides, where one makes contact with the area of higher permeability after 7 minutes. This contact causes fluid leakage and slows the fracture growth rate. Additionally, the tip in contact with the higher permeability region propagates much more slowly, while the remaining tip propagates in the direction of the maximum principal stress. In both cases studied, the fractures could not propagate through the higher permeability region. This example clearly illustrates the sensitivity of matrix permeability and demonstrates that simulation tools that do not capture this effect or represent flow in the matrix at all cannot accurately represent the propagation process.

### 4.3. Effect of injection rate

This example investigates the effect of the injection rate on the expansion of the stimulation area. Figure 9 illustrates a 2D fractured domain containing three fractures with boundary conditions described in the figure. We assume that the permeability is isotropic and homogeneous, given by $\kappa_{xx} = \kappa_{yy} = 5.0 \times 10^{-20}$ m$^2$. The principal stress is given by $\sigma_1 = 2\sigma_2 = 20$ MPa. Various injection rates are studied, and the effect on fracture growth and pressure in the fracture is shown in Figure 10.

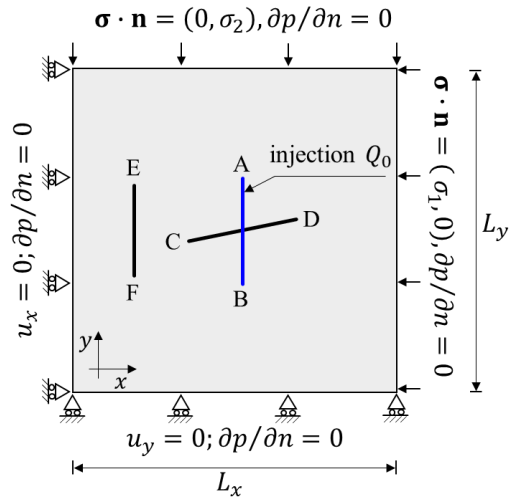

Figure 9. The geometry of model 3.

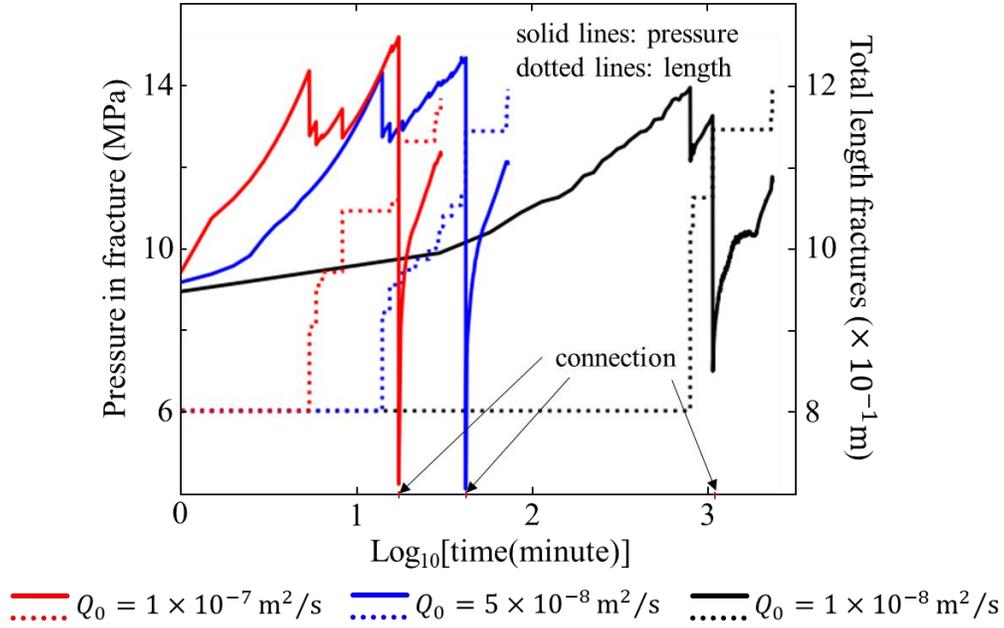

Figure 10. Effect of fluid injection rate on pressure at the injection point and total fracture growth.

The results shown in this example indicate that an increase in the injection rate leads to faster fracture propagation, and the propagation speed is nonlinearly dependent on the injection rate. As illustrated in Figure 10, wing cracks initiate after 4.5 minutes for an injection rate of $Q_0 = 1 \times 10^{-7}$ m²/s, whereas it takes up to 870 minutes for an injection rate of $Q_0 = 1 \times 10^{-8}$ m²/s. This indicates that increasing the injection rate by a factor of ten can accelerate the expansion of the fracture network by up to 200 times. However, if the injection rate is too low, then no fracture deformation may occur during our implementation.

## 4.4. Interaction with pre-existing fractures

Finally, we investigate the influence of the location and shape of pre-existing fractures on the expansion of the fracture network. The model geometry is shown in Figure 11. The matrix permeability in this example is assumed to be isotropic and homogeneous, i.e., $\kappa_{xx} = \kappa_{yy} = 5.0 \times 10^{-20}$ m². The principal stress is given by $\sigma_1 = 2\sigma_2 = 20$ MPa. The injection rate is $Q_0 = 2 \times 10^{-7}$ m²/s. The evolutions of the fracture geometry and the pore pressure are shown in Figure 13.

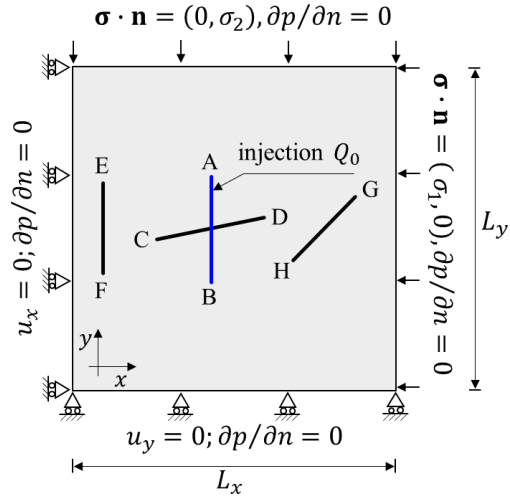

Figure 11. The geometry of model 4.

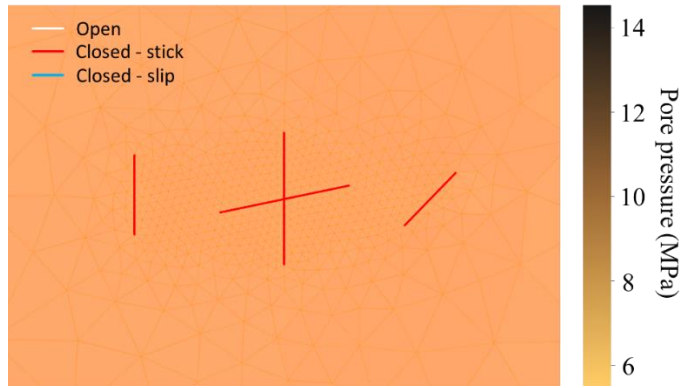

Figure 12. Fracture state and pressure in a 2D porous medium before fluid is injected.

Prior to fluid injection, the fracture mode and pore pressure are evaluated. As illustrated in Figure 12, pre-existing fractures are closed and remain in stick mode due to compressive stress and friction at the fracture interfaces. Additionally, the pressure throughout the domain is uniform at 6.8 MPa. The result in this simulation indicates a stable condition with no fracture slip or propagation.

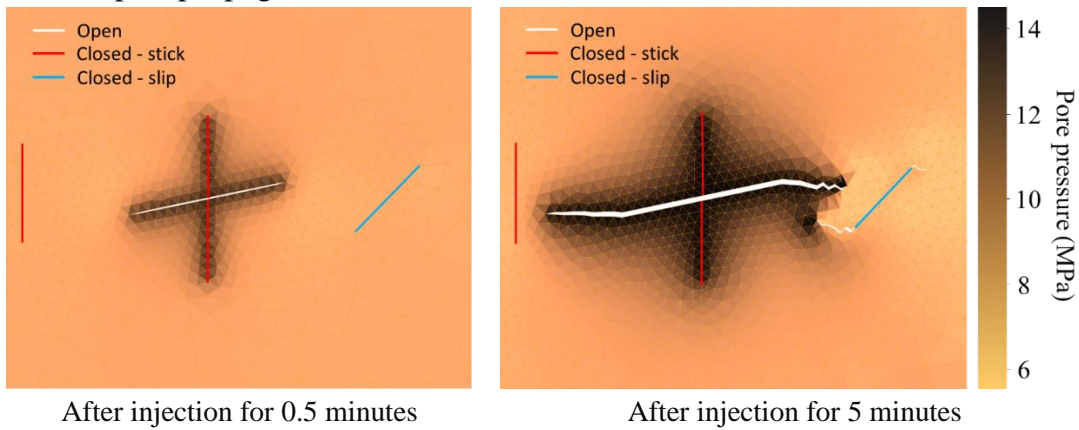

After injection for 0.5 minutes          After injection for 5 minutes

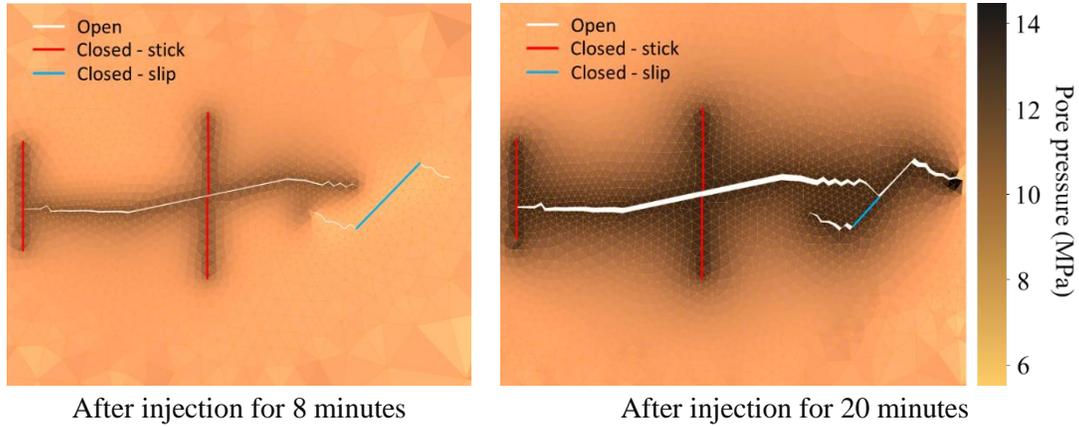

<div align="center">After injection for 8 minutes        After injection for 20 minutes</div>

Figure 13. Fracture propagation and pressure evolution in a 2D porous medium during fluid injection, $Q_0 = 2 \times 10^{-7} \, m^2/s$, into a pre-existing fracture. The solid white lines indicate open fractures, while the solid red lines indicate closed fractures. The color bar represents pore pressure in MPa.

Subsequently, fluid is injected into a vertical fracture, resulting in several interesting phenomena, as shown in Figure 13. First, the injection has an insignificant effect on the state of the fracture where fluid is injected, as it remains under compression under the influence of the stress regime. However, the injection facilitates the opening of the horizontal fracture connected to it and leads to the propagation of this fracture. Second, due to deformation and hydromechanical stress changes caused by fluid injection, the pre-existing fracture to the right of the domain starts to slip at an early stage of fluid injection. Eventually, small wing cracks are observed to form at the tips of this fracture. Third, there is a strong link between fracture propagation and pressure drop in the fracture. After a period of fluid injection, the pressure in the central, nearly horizontal, fracture increases sufficiently to cause tensile propagation of the fracture, which ultimately connects to the pre-existing fractures at the left and right. Each connection results in an instantaneous decrease in pore pressure, which takes time to recover through fluid injection before the fracture can resume growing. Furthermore, the expansion of the fractured network is influenced by the pre-existing fractures. During the simulation, the fracture on the right-hand side where the slip occurs continues to grow, while the fracture on the left side where compression occurs (closed in stick mode) prevents further network expansion.

## 5. Conclusions

This paper presents a mathematical model and numerical approach to investigate the use of mixed-mechanism stimulation to improve permeability in geothermal reservoirs. The mathematical model combines Biot poroelasticity and fracture mechanics and accounts for frictional contact mechanics and fracture propagation and connection. A two-level model that combines finite volume and finite element methods is proposed for numerical simulations. Several numerical examples are performed, and the results indicate the following:

1) Fluid injection at elevated pressure can induce shear slip and dilation, opening, and propagation of fractures. Newly formed fractures tend to propagate in the direction of maximum principal stress. In the case of multiple connected fractures in an anisotropic stress field, the propagation of fractures depends on fracture network characteristics such as fracture orientation relative to the stress field and whether fractures are hydraulically connected to the well through other fractures.

2) A more permeable bulk domain slows fracture growth by causing fluid leakage into the matrix, making hydraulic stimulations less effective for areas with higher permeability.

3) The relationship between the injection rate and fracture growth speed is nonlinear, and injection at a low rate may not result in fracture expansion. In most cases, when the injection rate is slower, the injection time required for a fracture to propagate is significantly longer.

4) The locations of pre-existing fractures influence the expansion of a fracture network. Fractures tend to propagate in the direction of the maximum principal stress, and pre-existing fractures can facilitate or impede the development of propagating fractures.

In conclusion, this study demonstrates that mixed-mechanism stimulation can significantly improve permeability by expanding the fracture network. However, this expansion is complex and influenced by various factors, including the stress state, material permeability, injection rate, and fracture location. The simulation model proposed in this study represents an approach that is appropriate for utilization in future studies to further investigate these phenomena.

## Declaration of Competing Interest

The authors declare that they have no known competing financial interests or personal relationships that could have appeared to influence the work reported in this paper.

## Acknowledgments

This project received funding from the European Research Council (ERC) under the European Union's Horizon 2020 research and innovation program (grant agreement No 101002507).